\documentclass[12pt]{article}
\newtheorem{theorem}{Theorem}[section]

\newtheorem{proposition}[theorem]{Proposition}

\newtheorem{corollary}[theorem]{Corollary}
\begin{document}
\begin{titlepage}

\title{Linear Connections on Light-like Manifolds}

%\author{\\
%\hspace{-1.4cm}\begin{tabular}{ccc} T.~Dereli\thanks{Electronic
%address:{\tt tdereli@ku.edu.tr}}&\c{S}.
%Ko\c{c}ak\thanks{Electronic address: {\tt
%skocak@anadolu.edu.tr}}&M. Limoncu\thanks{Electronic address:
%{\tt mlimoncu@anadolu.edu.tr}}\\
%\small{Department of Physics,}&\small{Department of
%Mathematics,}&\small{Department of Mathematics,}\\
%\small{Ko\c{c} University}&\small{Anadolu
%University}&\small{Anadolu University}\\
%\small{34450 Sar{\i}yer-\.{I}stanbul, Turkey}&\small{26470
%Eski\c{s}ehir, Turkey}&\small{26470 Eski\c{s}ehir, Turkey}
%\end{tabular}}

\author{ T. Dereli${}^{1}$\footnote{E.mail: tdereli@ku.edu.tr} , \c{S}.
Ko\c{c}ak${}^{2}$\footnote{E.mail: skocak@anadolu.edu.tr} , M.
Limoncu${}^{2}$\footnote{E.mail: mlimoncu@anadolu.edu.tr}\\
\\ ${}^{1}${\small Department of Physics, Ko\c{c} University, 34450
Sar{\i}yer-\.{I}stanbul, Turkey}\\
${}^{2}${\small Department of Mathematics, Anadolu University,
26470 Eski\c{s}ehir, Turkey}}

\maketitle

\begin{abstract}
\noindent It is well-known that a torsion-free linear connection
on a light-like manifold $(M,g)$ compatible with the degenerate
metric $g$ exists if and only if $Rad(TM)$ is a Killing
distribution. In case of existence, there is an infinitude of
connections with none distinguished. We propose a method to single
out connections with the help of a special set of 1-forms by the
condition that the 1-forms become parallel with respect to this
connection. Such sets of 1-forms could be regarded as an
additional structure imposed upon the light-like manifold. We
consider also connections with torsion and with non-metricity on
light-like manifolds.
\end{abstract}
\end{titlepage}

\section{Introduction}
In the following we will adopt the terminology of the book
Duggal-Bejancu \cite{Duggal-Bejancu}. A light-like manifold
$(M,g)$ is a smooth manifold $M$ with a smooth symmetric tensor
field of type (0,2) with constant index and nullity-degree
(co-rank) on $M$. There are half a dozen other names for such
manifolds, ``degenerate manifolds'' being maybe one of the most
popular terms. There are some scattered (and only partly related
and sometimes duplicated) works of mathematical and physical
origin in the literature about connections on light-like manifolds
(see \cite{Duggal-Bejancu} and references therein and also
\cite{Dombrowski-Horneffer}, \cite{Crampin},
\cite{Kunzle},\cite{Kozlov}).

It is an important result (it could be called the fundamental
theorem of connections on light-like manifolds) that a
torsion-free  linear connection $\nabla$ on $M$ compatible with
$g$ ($\nabla g=0$) exists if and only if $Rad(TM)$ is a Killing
distribution ([DB]). Here $Rad(TM)$ denotes the radical of $g$,
that is the sub-bundle of $TM$ with
$(Rad(TM))_{x}=Rad(T_{x}M)=\{\xi\in T_{x}M|\,g(\xi,v)=0\,\forall\,
v\in T_{x}M\}$ for $x\in M$. $Rad(TM)$ is a distribution of rank
equal to the constant nullity-degree (co-rank) of $g_{x}$. A
distribution $D$ on $M$ is called a Killing distribution if
$L_{X}g=0$ for each vector field $X\in\Gamma(D)$ ($\Gamma$ section
space), $L$ being the Lie derivative.

To motivate our approach we first consider a light-like manifold
$(M,g)$ with nullity degree 1. In this case , $Rad(TM)$ is a
1-dimensional distribution (a line bundle) on $M$. We assume this
line bundle to be trivial and choose a nowhere vanishing vector
field $\xi\in\Gamma(Rad(TM))$, i.e. a trivialization of this line
bundle. (Such a vector field is regarded by physicists to be a
``time-vector field''.) Additionally, we consider a 1-form field
$\tau\in\Gamma(T^{*}M)$ such that $\tau(\xi)=1$, i.e.
$\tau_{x}(\xi_{x})=1$ for all $x\in M$. It can easily be seen that
$\bar{g}:=g+\tau\otimes\tau$ is a non-degenerate symmetric
(0,2)-tensor field so that $(M,\bar{g})$ becomes a semi-Riemannian
manifold.
\begin{proposition}
i) Let $\nabla$ be a torsion-free linear connection on $M$ with
$\nabla_{X}g=0$ and $\nabla_{X}\tau=0$ for all $X\in\Gamma(TM)$.
Then, $L_{\xi}g=0$, $d\tau=0$ and $\nabla$ is the Levi-Civita
connection on $(M,\bar{g})$. ii) Assume $L_{\xi}g=0$ and
$d\tau=0$. If $\nabla$ is the Levi-Civita connection on
$(M,\bar{g})$, then $\nabla_{X}g=0$ and $\nabla_{X}\tau=0$ for all
$X\in\Gamma(TM)$.
\end{proposition}

We give below a proof of this proposition independently of the
above referred fundamental theorem. The following corollary
illustrates our viewpoint: Selection of a unique connection from
the infinitude of connections compatible with $g$ with the help of
an appropriate 1-form.
\begin{corollary}
Let $L_{\xi}g=0$ and let there exist a closed 1-form field $\tau$
with $\tau(\xi)=1$. Then there exists a unique torsion-free linear
connection $\nabla$ on $M$ with $\nabla_{X}g=0$ and
$\nabla_{X}\tau=0$ for all $X\in\Gamma(TM)$.
\end{corollary}
{\bf Proof} (of corollary) By part ii) of the proposition, at
least one connection of the desired properties exists: The
Levi-Civita of $\bar{g}=g+\tau\otimes\tau$. On the other hand, if
any torsion-free $\nabla$ satisfies $\nabla_{X}g=0$ and
$\nabla_{X}\tau=0$, then
$\nabla_{X}\bar{g}=\nabla_{X}g+\nabla_{X}\tau\otimes\tau+\tau\otimes\nabla_{X}\tau=0$,
so that $\nabla$ is necessarily the Levi-Civita of $\bar{g}$.
\newline
{\bf Proof} (of the proposition) i) It is obvious that $\Gamma$ is
the Levi-Civita connection on $(M,\bar{g})$ (as seen in the proof
of the corollary).\newline First we show $d\tau=0$:
\begin{eqnarray*}
(\nabla_{X}\tau)(Y)=\nabla_{X}(\tau(Y))-\tau(\nabla_{X}Y)=0\\
(\nabla_{Y}\tau)(X)=\nabla_{Y}(\tau(X))-\tau(\nabla_{Y}X)=0
\end{eqnarray*}
for $X,Y\in \Gamma(TM)$ by assumption. Subtracting these
equalities we get
\begin{eqnarray*}
X(\tau(Y))-Y(\tau(X))=\tau(\nabla_{X}Y-\nabla_{Y}X)=\tau([X,Y])
\end{eqnarray*}
which means $(d\tau)(X,Y)=0$.\newline Next we will show
$L_{\xi}g=0$, but we first remark that $\nabla_{X}\xi=0$ for
$X\in\Gamma(TM)$: $\bar{g}=g+\tau\otimes\tau$ gives
$\bar{g}(\xi,Y)=\tau(Y)$ for $Y\in\Gamma(TM)$. Then
$\nabla_{X}(\bar{g}(\xi,Y))=\nabla_{X}(\tau(Y))$,
\begin{eqnarray*}
(\nabla_{X}\bar{g})(\xi,Y)+\bar{g}(\nabla_{X}\xi,Y)+\bar{g}(\xi,\nabla_{X}Y)
=(\nabla_{X}\tau)(Y)+\tau(\nabla_{X}Y).
\end{eqnarray*}
$\nabla_{X}\bar{g}=0$, $\nabla_{X}\tau=0$ and
$\bar{g}(\xi,\nabla_{X}Y)=\tau(\nabla_{X}Y)$ imply
$\bar{g}(\nabla_{X}\xi,Y)=0$ which gives $\nabla_{X}\xi=0$ by
non-degeneracy of $\bar{g}$.\newline Now,
\begin{eqnarray*}
(\nabla_{X}g)(Y,Z)=\nabla_{X}(g(Y,Z))-g(\nabla_{X}Y,Z)-g(Y,\nabla_{X}Z)=0
\end{eqnarray*}
by assumption.
\begin{eqnarray*}
X(g(Y,Z))-g([X,Y]+\nabla_{Y}X,Z)-g(Y,[X,Z]+\nabla_{Z}X)=0\\
X(g(Y,Z))-g(L_{X}Y,Z)-g(Y,L_{X}Z)-g(\nabla_{Y}X,Z)-g(Y,\nabla_{Z}X)=0.
\end{eqnarray*}
If we now take $X=\xi$, the last two terms vanish and we get
\begin{eqnarray*}
\xi(g(Y,Z))-g(L_{\xi}Y,Z)-g(Y,L_{\xi}Z)=0
\end{eqnarray*}
which means $(L_{\xi}g)(Y,Z)=0$

ii) It will be enough to show $\nabla_{X}\tau=0$ as
$\nabla_{X}g=0$ is then a consequence of
$\nabla_{X}\bar{g}=0$.\newline We have
\begin{eqnarray*}
(\nabla_{X}\tau)(Y)=\nabla_{X}(\tau(Y))-\tau(\nabla_{X}Y)
\end{eqnarray*}
for $Y\in\Gamma(TM)$ and
\begin{eqnarray*}
(\nabla_{X}\tau)(Y)=X(\bar{g}(\xi,Y))-\bar{g}(\xi,\nabla_{X}Y))
\end{eqnarray*}
since $\bar{g}(\xi,Y)=\tau(Y)$ for any $Y$. From the Koszul
formula
\begin{eqnarray*}
2\bar{g}(\nabla_{X}Y,Z)&=&X\bar{g}(Y,Z)+Y\bar{g}(Z,X)-Z\bar{g}(X,Y)\\
&&-\bar{g}(X,[Y,Z])+\bar{g}(Y,[Z,X])+\bar{g}(Z,[X,Y])
\end{eqnarray*}
we can write
\begin{eqnarray*}
(\nabla_{X}\tau)(Y)&=&\frac{1}{2}(X\bar{g}(\xi,Y)-Y\bar{g}(\xi,X)-\bar{g}(\xi,[X,Y]))\\
&&+\frac{1}{2}(\xi\bar{g}(X,Y)-\bar{g}(X,[\xi,Y])-\bar{g}(Y,[\xi,X])).
\end{eqnarray*}
The first term vanishes because of $d\tau=0$:
\begin{eqnarray*}
(d\tau)(X,Y)&=&X\tau(Y)-Y\tau(X)-\tau([X,Y])=0\\
&=&X\bar{g}(\xi,Y)-Y\bar{g}(\xi,X)-\bar{g}(\xi,[X,Y])=0.
\end{eqnarray*}
For the second term, it will be enough to show $L_{\xi}\bar{g}=0$:
\begin{eqnarray*}
(L_{\xi}\bar{g})(X,Y)=\xi\bar{g}(X,Y)-\bar{g}(L_{\xi}X,Y)-\bar{g}(X,L_{\xi}Y).
\end{eqnarray*}
To show $L_{\xi}\bar{g}=0$, we first see $L_{\xi}\tau=0$:
\begin{eqnarray*}
(L_{\xi}\tau)(X)&=&L_{\xi}(\tau(X))-\tau(L_{\xi}X)\\
&=&\xi(\tau(X))-\tau([\xi,X])\\
&=&X(\tau(\xi))\\
&=&0
\end{eqnarray*}
by $d\tau=0$ and $\tau(\xi)=1$. Now,
$L_{\xi}\bar{g}=L_{\xi}g+(L_{\xi}\tau)\otimes\tau+\tau\otimes(L_{\xi}\tau)=0$.

This proposition can be extended in several directions. The
nullity degree (co-rank) of $g$ can be higher (multi-time models
in physics), the connection can have torsion and/or non-metricity
$Q$. We give below a theorem comprising all these cases.

Recall that a connection $\nabla$ is said to have non-metricity
$Q$, if $(\nabla_{Z}g)(X,Y)=Q(Z,X,Y)$ for $X,Y,Z\in \Gamma(TM)$ It
follows from the Koszul-method that for a given non-degenerate
metric on $M$, an anti-symmetric (1,2)-tensor $T(X,Y)$ and a
(0,3)-tensor $Q(Z,X,Y)$ symmetric in $X$ and $Y$, there exists a
unique linear connection on $TM$ having torsion $T$ and
non-metricity $Q$.
\begin{theorem}
Let $(M,g)$ be a light-like manifold with constant index and
nullity degree $r$ ($Rad(TM)$ is an $r$-distribution). Assume
$Rad(TM)$ to be a trivial bundle with everywhere independent
vector fields $\xi_{i}\in\Gamma(Rad(TM))$, $i=1,...,r$. Let
$\tau_{i}$ ($i=1,...,r$) be a set of 1-form fields satisfying
$\tau_{i}(\xi_{j})=\delta_{ij}$, $i,j=1,...,r$.
$\bar{g}:=g+\displaystyle\sum_{k=1}^{r}\tau_{k}\otimes\tau_{k}$ is
then a non-degenerate metric on $M$.

i) Let $\nabla$ be a linear connection on $M$ with torsion $T$
($T(X,Y)=\nabla_{X}Y-\nabla_{Y}X-[X,Y]$ for $X,Y\in\Gamma(TM)$)
and with non-metricity (0,3)-tensor $Q$, i.e.
$(\nabla_{Z}g)(X,Y)=Q(Z,X,Y)$. Assume additionally that the
1-forms $\tau_{i}$ are parallel with respect to $\nabla$:
$\nabla_{X}\tau_{i}=0$ for $X\in\Gamma(TM)$ and $i=1,...,r$.

Then the following relationships hold:
\begin{eqnarray}
(d\tau_{i})(X,Y)&=&\tau_{i}(T(X,Y))\label{qt}\\
(L_{\xi_{i}}g)(X,Y)&=&g(X,T(\xi_{i},Y))+g(Y,T(\xi_{i},X))+Q(\xi_{i},X,Y)\label{TQ}\\
&&-Q(X,Y,\xi_{i})-Q(Y,X,\xi_{i})\nonumber
\end{eqnarray}
for $X,Y\in\Gamma(TM)$, $i=1,...,r$.

ii) Let $T$ be an anti-symmetric (1,2)-tensor field and $Q$ a
(0,3)-tensor field (symmetric in the last two variables)
satisfying (\ref{qt}) and (\ref{TQ}).

If $\nabla$ is the connection on $(M,\bar{g})$ with torsion $T$
and non-metricity $Q$ (i.e. $(\nabla_{Z}\bar{g})(X,Y)=Q(Z,X,Y)$),
then $\nabla$ has also non-metricity $Q$ with respect to $g$ (i.e.
$(\nabla_{Z}g)(X,Y)=Q(Z,X,Y)$) and the 1-forms $\tau_{i}$ are
parallel with respect to $\nabla$ for $i=1,...,r$.
\end{theorem}
\begin{corollary}
In the setting of the above theorem, given any torsion tensor $T$
and non-metricity tensor $Q$ on $M$ satisfying (\ref{qt}) and
(\ref{TQ}), there exists a unique linear connection on $M$ having
torsion $T$, non-metricity $Q$ for $g$ and making $\tau_{i}$
($i=1,...,r$) parallel.
\end{corollary}

It might be useful to re-word the corollary for the
metric-compatible case also:
\begin{corollary}
In the setting of the above theorem, given any torsion tensor $T$
satisfying
\begin{eqnarray*}
(d\tau_{i})(X,Y)=\tau_{i}(T(X,Y))
\end{eqnarray*}
and
\begin{eqnarray*}
(L_{\xi_{i}}g)(X,Y)=g(X,T(\xi_{i},Y))+g(Y,T(\xi_{i},X))
\end{eqnarray*}
for $i=1,...,r$, there exists a unique linear connection on $M$
compatible with $g$, having torsion $T$ and making $\tau_{i}$
($i=1,...,r$) parallel.
\end{corollary}
As a last corollary, we note the torsion-free case also:
\begin{corollary}
In the setting of the above theorem, assume the 1-forms $\tau_{i}$
to be closed and $L_{\xi_{i}}g=0$ for $i=1,...,r$. Then there
exists a unique torsion-free and $g$-compatible connection
$\nabla$ with $\nabla\tau_{i}=0$ ($i=1,...,r$).
\end{corollary}
{\bf Proof} (of the theorem) i) From the equalities
\begin{eqnarray*}
(\nabla_{X}\tau_{i})(Y)=X\tau_{i}(Y)-\tau_{i}(\nabla_{X}Y)=0\\
(\nabla_{Y}\tau_{i})(X)=Y\tau_{i}(X)-\tau_{i}(\nabla_{Y}X)=0\\
\end{eqnarray*}
we get
\begin{eqnarray*}
X\tau_{i}(Y)-Y\tau_{i}(X)-\tau_{i}([X,Y]+T(X,Y))=0
\end{eqnarray*}
which means
\begin{eqnarray*}
(d\tau_{i})(X,Y)=X\tau_{i}(Y)-Y\tau_{i}(X)-\tau_{i}([X,Y])=\tau_{i}(T(X,Y))
\end{eqnarray*}
so that the first relationship (\ref{qt}) holds.\newline Now, let
us start with
\begin{eqnarray}\label{v}
(\nabla_{Z}g)(X,Y)=Zg(X,Y)-g(\nabla_{Z}X,Y)-g(X,\nabla_{Z}Y)=Q(Z,X,Y)
\end{eqnarray}
and insert $\xi_{i}$ for $Z$ (using torsion):
\begin{eqnarray*}
\xi_{i}g(X,Y)-g(\nabla_{X}\xi_{i}+[\xi_{i},X]+T(\xi_{i},X),Y)\hspace{-5cm}\\
&-g(X,\nabla_{Y}\xi_{i}+[\xi_{i},Y]+T(\xi_{i},Y))=Q(\xi_{i},X,Y),
\end{eqnarray*}
\begin{eqnarray*}
\xi_{i}g(X,Y)-g([\xi_{i},X],Y)
-g(X,[\xi_{i},Y])=Q(\xi_{i},X,Y)+g(\nabla_{X}\xi_{i},Y)+g(X,\nabla_{Y}\xi_{i})\hspace{-7.3cm}\\
&+g(T(\xi_{i},X),Y)+g(X,T(\xi_{i},Y)).
\end{eqnarray*}
As $(L_{\xi_{i}}g)(X,Y)$ equals the left-hand side, it is enough
to see $$g(\nabla_{X}\xi_{i},Y)=-Q(X,Y,\xi_{i})$$ for any
$X,Y\in\Gamma(TM)$ for the second relationship to hold. Now,
inserting $\xi_{i}$ for $Y$ in (\ref{v}), we get
\begin{eqnarray*}
Zg(X,\xi_{i})-g(\nabla_{Z}X,\xi_{i})-g(X,\nabla_{Z}\xi_{i})=Q(Z,X,\xi_{i})
\end{eqnarray*}
which reduces to $-g(X,\nabla_{Z}\xi_{i})=Q(Z,X,\xi_{i})$ for any
$X,Z\in\Gamma(TM)$ since $\xi_{i}\in\Gamma(Rad(TM))$.

ii)We first note $\bar{g}(\xi_{i},Y)=\tau_{i}(Y)$ for
$Y\in\Gamma(TM)$ by definition of $\bar{g}$.

The following equations can be verified by direct computation
using (\ref{qt}) and (\ref{TQ}):
\begin{eqnarray*}
\hspace{1.2cm}X\tau_{i}(Y)-\tau_{i}(\nabla_{X}Y)=\frac{1}{2}(X\bar{g}(\xi_{i},Y)-Y\bar{g}(\xi_{i},X)+\xi_{i}\bar{g}(X,Y)
\hspace{-8.1cm}\\
&+\bar{g}(X,[Y,\xi_{i}])-\bar{g}(Y,[\xi_{i},X])-\bar{g}(\xi_{i},[X,Y])\hspace{-0.3cm}\\
&+Q(X,Y,\xi_{i})+Q(Y,X,\xi_{i})-Q(\xi_{i},X,Y)\\
&\qquad\qquad+\bar{g}(X,T(Y,\xi_{i}))+\bar{g}(Y,T(X,\xi_{i}))-\bar{g}(\xi_{i},T(X,Y)))
\end{eqnarray*}
\begin{eqnarray*}
\hspace{4.82cm}=\frac{1}{2}(X\tau_{i}(Y)-Y\tau_{i}(X)-\tau_{i}([X,Y])
\hspace{-7.4cm}\\
&+\xi_{i}\bar{g}(X,Y)-\bar{g}([\xi_{i},X],Y)-\bar{g}(X,[\xi_{i},Y])\hspace{-0.1cm}\\
&+Q(X,Y,\xi_{i})+Q(Y,X,\xi_{i})-Q(\xi_{i},X,Y)\\
&\qquad\qquad+\bar{g}(X,T(Y,\xi_{i}))+\bar{g}(Y,T(X,\xi_{i}))-\tau_{i}(T(X,Y)))
\end{eqnarray*}
\begin{eqnarray*}
\hspace{-0.5cm}X\tau_{i}(Y)-\tau_{i}(\nabla_{X}Y)=\frac{1}{2}((d\tau_{i})(X,Y)-\tau_{i}(T(X,Y))+(L_{\xi_{i}}\bar{g})(X,Y)
\hspace{-9.3cm}\\
&+Q(X,Y,\xi_{i})+Q(Y,X,\xi_{i})-Q(\xi_{i},X,Y)\hspace{-2.8cm}\\
&\qquad\qquad+\bar{g}(X,T(Y,\xi_{i}))+\bar{g}(Y,T(X,\xi_{i}))).
\end{eqnarray*}
Thus we get
\begin{eqnarray}
\hspace{0.5cm}(\nabla_{X}\tau_{i})(Y)=\frac{1}{2}((L_{\xi_{i}}\bar{g})(X,Y)
+Q(X,Y,\xi_{i})+Q(Y,X,\xi_{i})\hspace{-7.5cm}\label{w}\\
&-Q(\xi_{i},X,Y)+\bar{g}(X,T(Y,\xi_{i}))+\bar{g}(Y,T(X,\xi_{i}))).\nonumber
\end{eqnarray}
Let us now compute $L_{\xi_{i}}\bar{g}$ using
$\bar{g}=g+\displaystyle\sum_{k=1}^{r}\tau_{k}\otimes\tau_{k}$:
\begin{eqnarray*}
(L_{\xi_{i}}\bar{g})(X,Y)=(L_{\xi_{i}}g)(X,Y)+\sum_{k=1}^{r}((L_{\xi_{i}}\tau_{k})(X)\tau_{k}(Y)
+\tau_{k}(X)(L_{\xi_{i}}\tau_{k})(Y))
\end{eqnarray*}
\begin{eqnarray*}
(L_{\xi_{i}}\tau_{k})(X)&=&\xi_{i}\tau_{k}(X)-\tau_{k}([\xi_{i},X])\\
&=&\xi_{i}\tau_{k}(X)-X\tau_{k}(\xi_{i})-\tau_{k}([\xi_{i},X])\\
&=&(d\tau_{k})(\xi_{i},X)\\
&=&\tau_{k}(T(\xi_{i},X))
\end{eqnarray*}
since $\tau_{k}(\xi_{i})=\delta_{ki}$ and the last equality by
assumption (\ref{qt}). Thus we obtain
\begin{eqnarray*}
(L_{\xi_{i}}\bar{g})(X,Y)&=&(L_{\xi_{i}}g)(X,Y)+\sum_{k=1}^{r}(\tau_{k}(T(\xi_{i},X))\tau_{k}(Y)
+\tau_{k}(X)\tau_{k}(T(\xi_{i},Y)))\\
&=&(L_{\xi_{i}}g)(X,Y)-\sum_{k=1}^{r}\tau_{k}(X)\tau_{k}(T(Y,\xi_{i}))
-\sum_{k=1}^{r}\tau_{k}(Y)\tau_{k}(T(X,\xi_{i}))\\
&=&(L_{\xi_{i}}g)(X,Y)-\sum_{k=1}^{r}(\tau_{k}\otimes\tau_{k})(X,T(Y,\xi_{i}))
-\sum_{k=1}^{r}(\tau_{k}\otimes\tau_{k})(Y,T(X,\xi_{i}))
\end{eqnarray*}
\begin{eqnarray*}
\hspace{1.6cm}&=&-Q(X,Y,\xi_{i})-Q(Y,X,\xi_{i})+Q(\xi_{i},X,Y)\\
&&-g(X,T(Y,\xi_{i}))-g(Y,T(X,\xi_{i}))\\
&&-\sum_{k=1}^{r}(\tau_{k}\otimes\tau_{k})(X,T(Y,\xi_{i}))
-\sum_{k=1}^{r}(\tau_{k}\otimes\tau_{k})(Y,T(X,\xi_{i}))\\
&=&-Q(X,Y,\xi_{i})-Q(Y,X,\xi_{i})+Q(\xi_{i},X,Y)\\
&&-\bar{g}(X,T(Y,\xi_{i}))-\bar{g}(Y,T(X,\xi_{i})).
\end{eqnarray*}
Inserting this into the equation (\ref{w}) we get
\begin{eqnarray*}
(\nabla_{X}\tau_{i})(Y)=0
\end{eqnarray*}
The other assertion, $(\nabla_{Z}g)(X,Y)=Q(Z,X,Y)$ is now a
consequence of
\begin{eqnarray*}
\nabla_{Z}\bar{g}=\nabla_{Z}g+\nabla_{Z}(\sum_{k=1}^{r}\tau_{k}\otimes\tau_{k})=\nabla_{Z}g.
\end{eqnarray*}
\newline

\end{document}